\documentclass[11pt]{article}
\usepackage{amssymb}
\usepackage{amsmath}

\newtheorem{theorem}{Theorem}
\newtheorem{lemma}{Lemma}

\begin{document}

\title{Local Existence for the One-dimensional Vlasov-Poisson System with Infinite Mass}
\author{Stephen Pankavich \\
Department of Mathematics \\
Indiana University \\
Bloomington, IN 47401 \\
sdp@indiana.edu}
\date{\today}
\maketitle

\begin{center}
\end{center}

\begin{abstract}
A collisionless plasma is modelled by the Vlasov-Poisson system in
one dimension.  We consider the situation in which mobile negative
ions balance a fixed background of positive charge, which is
independent of space and time, as $\vert x \vert \rightarrow
\infty$.  Thus, the total positive charge and the total negative
charge are both infinite.  Smooth solutions with appropriate
asymptotic behavior are shown to exist locally in time, and
criteria for the continuation of these solutions are established.
\end{abstract}

\section*{Introduction}

Consider the following simplified model of a collisionless plasma,
using the Vlasov-Poisson system in one dimension. Let $F :
\mathbb{R} \rightarrow [0,\infty)$ and $f_0 : \mathbb{R} \times
\mathbb{R} \rightarrow [0,\infty)$ be given, and seek a function $f
: [0,\infty) \times \mathbb{R} \times \mathbb{R} \rightarrow
[0,\infty)$ such that

\begin{equation}
\label{VP} \left. \begin{array}{ccc}
& & \partial_t f + v \ \partial_x f - E \ \partial_v f = 0\\
\\
& & \rho(t,x) = \int \left ( F(v) - f(t,x,v) \right ) \ dv \\
\\
& & E(t,x) = \frac{1}{2} \left ( \int_{-\infty}^x \rho(t,y) \ dy - \int_x^\infty \rho(t,y) \ dv \right )\\
\\
& & f(0,x,v) = f_0(x,v) \end{array} \right \}
\end{equation}

\addvspace{0.5in}

\noindent where $t \in [0,\infty)$ denotes time, $x \in \mathbb{R}$
denotes space, and $v \in \mathbb{R}$ denotes momentum. Here, $f$
gives the density in phase space of mobile negative ions, while $F$
describes a number density of positive ions which form a fixed
background.  We seek solutions for which $f(t,x,v) \rightarrow F(v)$
as $\vert x \vert \rightarrow \infty$. Precise conditions which
ensure local-in-time existence of such solutions to (\ref{VP}) are
given in Section $1$.

To this end, we wish to proceed in a manner similar to \cite{VPSSA},
which proved the existence of a local-in-time solution to the three
dimensional analogue of (\ref{VP}).  In that paper, the existence
argument hinged upon showing that $\rho$ decayed faster than $\vert
x \vert^{-3}$.  In the same vein, the main difficulty of this paper
arises in showing $\rho$ decays rapidly enough in $\vert x \vert$,
in this case, better than $\vert x \vert^{-1}$. A crucial difference
between these arguments, however, lies in the decay rates of the
other functions. Let $g = F - f$.  In the three dimensional problem
it is well known that $E$, $g$ and $\nabla g$ decay at a rate of
$\vert x \vert^{-2}$. For (\ref{VP}), none of the analogous terms
are initially known to tend to zero for large $x$. Indeed, of these
terms only $\partial_x g$ can be shown to do so. One important part
of the existence proof in \cite{VPSSA} relied upon Lemma $1$ and $5$
of that paper generating the estimate :
$$ \rho \approx \int E \cdot \nabla g \ dv \approx \vert x \vert^{-4}.$$  Due
to the differing behavior of $E$ and $\nabla g$ in (\ref{VP}), this
estimate does not hold in the one dimensional problem. Instead,
since neither $E$ nor $\nabla g$ are known to decay in space, it is
unclear as to the a priori behavior of $\rho$. This difficulty is
remedied through the use of Lemma $4$. Regardless, it is worth
noting that due to a loss of spatial decay in each function,
(\ref{VP}) contains additional difficulties to those encountered in
the three dimensional problem, and many of the
techniques used in \cite{VPSSA} cannot be utilized here.\\

The Vlasov-Poisson system has been studied extensively in the case
where $F(v) = 0$ and solutions tend to zero as $\vert x \vert
\rightarrow \infty$, mostly with respect to the three-dimensional
problem.  Most of the literature involving the one-dimensional
Vlasov-Poisson system focus on time asymptotics, such as \cite{BKR}
and \cite{BFFM}.  Much more work has been done concerning the
three-dimensional problem.  Smooth solutions were shown to exist
globally-in-time in \cite{Pfaf} (refined in \cite{Sch}) and
independently in \cite{LP}. Important results preliminary to the
discovery of a global-in-time solution include \cite{Batt} and
\cite{Horst}.  A complete discussion of the literature concerning
the Vlasov-Poisson system may be found in \cite{Glassey}, and more
recently \cite{Rein}.\\

A good deal of progress has also been made on the infinite mass
problem for the Vlasov-Poisson system, this being the
three-dimensional analogue of (\ref{VP}).  Other than \cite{VPSSA},
smooth solutions were shown to exist globally-in-time for the case
of a radial field in \cite{PankavichR}. A priori bounds on the
charge and current densities were established in \cite{SSAVP},
assuming $F(v)$ is suitably smooth, radial and decreasing. Finally,
using the bounds derived in \cite{SSAVP}, a unique global-in-time
smooth solution was shown to exist in \cite{PankavichG} for an
arbitrary electric field, assuming $F$ is compactly supported.
Hence, our assumptions on data will mimic that of \cite{VPSSA} and
\cite{PankavichG}.

\section*{Section 1}

We let $p > 1$ be given and use the following notation.  Denote
$$ R(x) = R(\vert x \vert) = \sqrt{1 + \vert x \vert^2}.$$
For given functions $h : \mathbb{R}^2 \rightarrow \mathbb{R}$ and
$\sigma : \mathbb{R} \rightarrow \mathbb{R}$, we will use the norms
$$ \Vert h \Vert_\infty = \sup_{z \in \mathbb{R}^2} \vert h(z) \vert,$$
$$\Vert \sigma \Vert_p = \Vert \sigma R^p(x) \Vert_{L^\infty(\mathbb{R})},$$
and
$$ \Vert \vert h \vert \Vert = \Vert h \Vert_\infty +
\Vert \partial_v h \Vert_\infty + \Vert
\partial_x h \Vert_p + \Vert \int h \ dv \Vert_p.$$
 but never use the $L^p$ norm (for $p$ finite).  For example,
we will write $\Vert \rho(t) \Vert_p$ for the $\Vert \cdot
\Vert_p$ norm of $x \rightarrow \rho(t,x)$.\\

In order to study the above system, we will assume the following
conditions hold throughout :
\begin{enumerate}
\renewcommand{\labelenumi}{(A-\arabic{enumi})}
\item $f_0 \in \mathcal{C}^1(\mathbb{R}^2)$ is nonnegative with
compact $v$-support and for $x,v \in \mathbb{R}$

$$\vert (F - f_0)(x,v) \vert  + \vert \partial_v (F - f_0)(x,v) \vert \leq CR^{-p}(x).$$

\item $F \in \mathcal{C}^3(\mathbb{R})$ is nonnegative and there
is $W \in (0,\infty)$ such that for $\vert v \vert > W$,
$$ F(v) = 0.$$

\item $\Vert \vert (F - f_0) \vert \Vert$ is finite.

\end{enumerate}
From these assumptions, local-in-time existence follows.

\begin{theorem}
Assuming (A-$1$) thru (A-$3$) hold, there exist $\delta > 0$ and
$f \in C^1([0,\delta] \times \mathbb{R}^2)$ satisfying (\ref{VP})
with $\Vert \vert (F - f)(t) \vert \Vert\leq C$ for $t \in
[0,\delta]$. Moreover, $f$ is unique.
\end{theorem}
\vspace{0.25in} In addition, we may continue the local-in-time
solution as long as this norm remains bounded.
\begin{theorem}
Assume (A-$1$) thru (A-$3$) hold.  Let $T > 0$ be given and $f$ be a
$C^1$ solution of (\ref{VP}) on $[0,T] \times \mathbb{R}^2$. If
$$ \Vert \vert (F-f)(t) \vert \Vert \leq C$$ on $[0,T]$,
then we may extend the solution to $[0,T + \delta]$ for some $\delta
> 0$ where $\Vert \vert (F-f)(t) \vert \Vert$ is bounded on $[0,T +
\delta]$.
\end{theorem}

\vspace{0.25in}

In order to combine the proofs of Theorems $1$ and $2$, let
$f^{(T)} = f_0$ if $T = 0$ and for $T > 0$, let $f^{(T)}$ be the
solution assumed to exist in Theorem $2$.  For $\delta > 0 $ and
$r > 1$, define

$$ \begin{array}{ccl} \mathcal{C} = \mathcal{C}(\delta, r) & := & \{ f \in C^1([0,T+\delta] \times
\mathbb{R}^2) : \\ \ & \ & \ f(t,x,v) = f^{(T)}(t,x,v) \ \mathrm{if}
\ t \in [0,T] \\ & \ & \ \mathrm{and} \ \Vert \vert (F - f)(t) \vert
\Vert \leq r \ \mathrm{if} \ t \in [0,T + \delta] \}.
\end{array}$$ We take $$r \geq 1 + \sup_{t \in [0,T]} \Vert \vert
(F - f)(t) \vert \Vert $$ and $0 \leq \delta \leq \frac{1}{r}$ for
the remainder of the paper.\\

For $f \in \mathcal{C}$, define \begin{eqnarray}
\label{g} g & = & F - f, \\
\label{rho} \rho & = & \int g \ dv, \\
\label{E} E & = & \frac{1}{2} \left ( \int_{-\infty}^x \rho(t,y) \
dy - \int_x^\infty \rho(t,y) \ dy \right ).
\end{eqnarray}
Further define $\mathcal{F}[f] = \tilde{f}$ by
\begin{equation}
\label{ftildedef}
\left \{
\begin{array}{rcl}
\tilde{f}(t,x,v) = f^{(T)}(t,x,v), & & t \in [0,T] \\
\partial_t \tilde{f} + v \ \partial_x \tilde{f} - E \ \partial_v \tilde{f} =
0, & & t \in [T,T  + \delta]
\end{array} \right .
\end{equation}
Then, $$ \tilde{g} = F - \tilde{f},$$ and $$\tilde{\rho} = \int
\tilde{g} \ dv.$$

In Section $2$, we will choose $r$ and $\delta$ such that
$\mathcal{F} : \mathcal{C} \rightarrow \mathcal{C}$, and use this
to show that an iterative sequence converges, thus proving
Theorems $1$ and $2$.  Unless it is stated otherwise, we will
denote by ``C" a generic constant which changes from line to line
and may depend upon $f_0, F, T$, or $\sup_{t \in [0,T]} \Vert
\vert (F - f)(t) \vert \Vert$, but not on $t,x,v,r$ or $\delta$.
When it is necessary to refer to specific constants, we will use
superscripts.  For
instance, $C^{(1)}$ will always denote the same value.\\
Define
$$ I(t - T) = \left \{ \begin{array}{ll} 0, & t \leq T \\ 1, & t > T. \end{array}
\right. $$ Then, by definition, $f \in \mathcal{C}$ implies
\begin{equation}
\label{rhobound} \Vert \rho(t) \Vert_p \leq C (1 + r I(t-T))
\end{equation}
for all $x \in \mathbb{R}$ and $t \in [0,T+\delta]$.  Therefore,
given $t \in [0,T + \delta]$ and $x \in \mathbb{R}$, we find
\begin{equation}
\label{Ebound} \vert E(t,x) \vert \leq \Vert \rho(t) \Vert_p \left
( \int R^{-p}(y) \ dy \right ) \leq C \Vert \rho(t) \Vert_p
. \\
\end{equation}
Thus, for $t \in [0,T+\delta]$, we find
\begin{equation}
\label{intE}
\int_0^t \Vert E(s) \Vert_\infty \ ds \leq C \int_0^t (
1 + r I(s-T) ) \ ds \leq C (1 + r \delta) \leq C =: C^{(1)}
\end{equation}
and the field integral is uniformly bounded. Define the
characteristics, $X(s,t,x,v)$ and $V(s,t,x,v)$, as solutions to the
system of ordinary differential equations :
\begin{equation}
\label{char} \left. \begin{array}{ccc} & &
\frac{\partial}{\partial s} X(s,t,x,v) = V(s,t,x,v) \\
& & \frac{\partial}{\partial s} V(s,t,x,v) = - E(s,X(s,t,x,v)) \\
& & X(t,t,x,v) = x \\
& & V(t,t,x,v) = v. \end{array} \right \}
\end{equation}

We will make use of the following lemma which results merely from
this bound on the time integral of the field.

\begin{lemma} Let $t \in [0,T + \delta]$, $s \in [0,t]$, and $x \in
\mathbb{R}$ be given.  Then, for any $v \in \mathbb{R}$, $$ \vert
v \vert - C^{(1)} \leq \vert V(s,t,x,v) \vert \leq \vert v \vert +
C^{(1)}.$$ In particular, for $\vert v \vert > 2 C^{(1)}$,
$$ \frac{1}{2} \vert v \vert \leq \vert V(s,t,x,v) \vert \leq
\frac{3}{2} \vert v \vert.$$

\end{lemma}

\vspace{0.5in}

The following lemma is the crucial tool used in showing decay (at
a rate of $\vert x \vert^{-p}$) of the charge density, $\rho$.

\begin{lemma}
Assume that $\mathcal{E} : \mathbb{R} \rightarrow \mathbb{R}$ is
$\mathcal{C}^1$ and there exists $B > 0$ such that
$$ \vert \mathcal{E}(x) \vert \leq B$$
and $$ \vert \mathcal{E}^\prime(x) \vert \leq BR^{-p}(x)$$ for all
$x \in \mathbb{R}$.  Also, assume $H: \mathbb{R} \rightarrow
\mathbb{R}$ is $\mathcal{C}_c^2$.  Then,
$$\left \vert \int \mathcal{E}(X(s,t,x,v))
H^\prime(V(s,t,x,v)) \ dv \right \vert \leq CBR^{-p}(x)$$ for all $x
\in \mathbb{R}$ and $0 \leq s \leq t \leq T + \delta$.
\end{lemma}
We will postpone the proofs of all lemmas until Section $3$.

\section*{Section 2}

Let us estimate $\Vert \tilde{\rho}(t) \Vert_p$.  Define
\begin{equation}
\label{D}
D = \partial_t + v \partial_x - E \partial_v.
\end{equation}
Then,
\begin{equation}
\label{Dg}
D\tilde{g} = - E F^\prime
\end{equation}
so that
\begin{equation}
\label{rhotilde}
\tilde{\rho}(t,x) = \rho_0(t,x) - \int_0^t \int
E(s,X(s,t,x,v)) F^\prime(V(s,t,x,v)) \ dv \ ds.
\end{equation} where
$$ \rho_0(t,x) = \int \left ( F - f_0 \right ) (X(0,t,x,v),
V(0,t,x,v)) \ dv.$$ Since $\partial_x E = \rho$, we use
(\ref{rhobound}) and (\ref{Ebound}) in Lemma $2$ to
find

\begin{equation}
\begin{array}{rcl}
\label{EFint} \left \vert \int_0^t \int E(s, X(s,t,x,v))
F^\prime(V(s,t,x,v)) \ dv \ ds \right \vert & \leq &
\int_0^t C \left( 1 + r I(s-T) \right ) R^{-p}(x) \ ds \\
\\
& \leq & C \left ( 1 + r \delta \right ) R^{-p}(x) \\
\\
& \leq & CR^{-p}(x).
\end{array}
\end{equation}

Then, to estimate $\rho_0$, we use (A-1) and (A-2) to conclude that
$F - f_0$ has compact $v$-support.  Let $P_V := \sup \{\vert v \vert
: \exists x \in \mathbb{R} \ \hbox{such that} \ (F- f_0)(x,v) \neq 0
\}$ and define $P_v := P_V + C^{(1)}$. Using Lemma $1$ for $\vert v
\vert > P_v$, we find $\vert V(0) \vert
> P_V$ and thus $(F-f_0)(X(0),V(0)) = 0$.  It follows that
$(F-f_0)(X(0),V(0))$ has compact support as a function of $v$.  Now,
assume $\vert v \vert \leq P_v$. For $\vert x \vert > 2(P_v +
C^{(1)}) (T + \delta)$, we again use Lemma $1$ to find
$$\vert X(0,t,x,v) \vert \geq \vert x \vert - \int_0^t \vert V(\tau)
\vert \ d\tau \geq \vert x \vert - (P_v + C^{(1)})(T + \delta) \geq
\frac{1}{2} \vert x \vert$$ and thus $R^{-p}(X(0)) \leq
R^{-p}(\frac{1}{2}\vert x \vert) \leq C R^{-p}(x)$. In addition, for
$\vert x \vert \leq 2(P_v + C^{(1)})(T + \delta)$, we have
$R^{-p}(X(0)) \leq C \leq CR^{-p}(x).$ Therefore, we find for any $x
\in \mathbb{R}$,
$$ R^{-p}(X(0)) \leq C R^{-p}(x).$$
Hence, we find
\begin{eqnarray*}
\vert \rho_0(t,x) \vert & \leq & \int_{\vert v \vert \leq P_v} C
R^{-p}(X(0,t,x,v)) \ dv \\
& \leq & C R^{-p}(x).
\end{eqnarray*}
Using this and (\ref{EFint}), we find
\begin{equation}
\label{rhotildebound} \vert \tilde{\rho}(t,x) \vert \leq
CR^{-p}(x)
\end{equation}
and thus
$$ \Vert \tilde{\rho}(t) \Vert_p \leq C.$$

Now, we estimate $\tilde{g}$, $\partial_v \tilde{g}$, and
$\partial_x \tilde{g}$. First, we use (A-2), (A-3), and
(\ref{intE}) in (\ref{Dg}), after integrating along
characteristics, to  find
\begin{eqnarray*}
\vert \tilde{g}(t,x,v) \vert & \leq & \vert (F - f_0)(X(0), V(0))
\vert +
\int_0^t \vert E(s,X(s)) \vert \ \vert F^\prime(V(s)) \vert \ ds \\
& \leq & C + \Vert F^\prime \Vert_\infty \left ( \int_0^t \vert
E(s, X(s)) \vert \ ds \right ) \\
& \leq & C.
\end{eqnarray*}

Using $D$ defined by (\ref{D}), we have
$$ D(\partial_x \tilde{g}) = \partial_x E \left ( \partial_v \tilde{g} -
F^\prime \right ) = \rho \left ( \partial_v \tilde{g} - F^\prime
\right )$$ and
\begin{equation}
\label{Ddvg}
D(\partial_v \tilde{g}) = - \left (
\partial_x g + E F^{\prime \prime} \right ).
\end{equation}
Then, using (\ref{rhobound}), (\ref{Ebound}), and (A-2), we find
\begin{equation}
\label{Dxg} \vert D(\partial_x \tilde{g}) \vert \leq C R^{-p}(x) (
1 + r I(t-T))(1 + \vert \partial_v \tilde{g} \vert )
\end{equation}
and \begin{equation} \label{Dvg} \vert D(\partial_v \tilde{g})
\vert \leq C ( 1 + r I(t-T))(1 + \vert \partial_x \tilde{g} \vert)
\end{equation}
Denoting $$\tilde{G}(s) = \left ( \vert \partial_v \tilde{g} \vert +
\vert \partial_x \tilde{g} \vert \right ) (s,X(s),V(s)),$$ we
combine (\ref{Dxg}) and (\ref{Dvg}) and integrate along
characteristics to find
$$\tilde{G}(t) \leq \tilde{G}(0) + C \int_0^t (1 + r I(s-T)) \ (
1+ \tilde{G}(s) ) \ ds.$$  By Gronwall's inequality and (A-3), we
find $$\tilde{G}(t) \leq C (1 + \tilde{G}(0)) \leq C,$$ and it
follows that for any $t \in [0,t+ \delta]$,
\begin{equation}
\label{gtildebound} \Vert \tilde{g}(t) \Vert_\infty + \Vert
\partial_v \tilde{g}(t) \Vert_\infty + \Vert
\partial_x \tilde{g}(t) \Vert_\infty \leq C.
\end{equation}

For the remainder of the section, we will need to estimate terms
which involve $v$-integrals of $\tilde{g}$.  Thus, the following
lemma will be useful :

\begin{lemma}
For any $f \in \mathcal{C}$, define $\tilde{f}$ by (\ref{ftildedef})
and $\tilde{g} = F- \tilde{f}$.  Then, $Q_{\tilde{g}}(t)$ defined by

\begin{equation} \label{Qg} Q_{\tilde{g}}(t) = \sup \{ \vert v \vert :
\exists \ x \in \mathbb{R}, \tau \in [0,t] \ \mathrm{such \ that} \
\tilde{g}(\tau,x,v) \neq 0 \}.
\end{equation}
is bounded for all $t \in [0,T + \delta]$.
\end{lemma}

\addvspace{0.5in}

Continuing with the estimate of $\Vert \vert \tilde{g} \vert \Vert$,
we use (\ref{Dxg}) and (\ref{gtildebound}), and integrate along
characteristics to find $$ \vert \partial_x \tilde{g}(t,x,v) \vert
\leq \vert
\partial_x \tilde{g}(0,X(0),V(0) \vert + C\int_0^t (1 + r I(s-T))
R^{-p}(X(s)) \ ds.$$ By the definition of $Q_{\tilde{g}}(t)$, we
find that for $\vert x \vert > 2(T + \delta) Q_{\tilde{g}}(T +
\delta)$, we have $$\vert X(s,t,x,v) \vert \geq \vert x \vert -
\int_s^t \vert V(\tau) \vert \ d\tau \geq \vert x \vert - (T +
\delta)Q_{\tilde{g}}(T + \delta) \geq \frac{1}{2} \vert x \vert$$
and thus $R^{-p}(X(s)) \leq R^{-p}(\frac{1}{2}\vert x \vert) \leq C
R^{-p}(x)$. In addition, for $\vert x \vert \leq 2(T +
\delta)Q_g(T+\delta)$, we have $R^{-p}(X(s)) \leq C \leq
CR^{-p}(x).$ Therefore, we find for any $x \in \mathbb{R}$,
$$ R^{-p}(X(s)) \leq C R^{-p}(x).$$ We use this and (A-3) to find
$$ \vert \partial_x \tilde{g}(t,x,v) \vert \leq C R^{-p}(x),$$
and therefore
\begin{equation}
\label{Dxgbound} \Vert \partial_x \tilde{g}(t,x,v) \Vert_p \leq C.
\end{equation}
Collecting (\ref{rhotildebound}), (\ref{gtildebound}), and
(\ref{Dxgbound}) we find
\begin{equation}
\label{gtildetriple} \Vert \vert \tilde{g}(t,x,v) \vert \Vert \leq
C^{(2)}.
\end{equation}
We take $r \geq C^{(2)}$ and conclude that $\tilde{f} \in
\mathcal{C}$ and $\mathcal{F} : \mathcal{C}
\rightarrow \mathcal{C}.$\\

Next, let $f,\rho,E, \tilde{f}$, and $\tilde{\rho}$ be as before and
for $h \in \mathcal{C}$, define
$$ g_h = F - h,$$
$$\rho_h = \int g_h \ dv$$
$$ E_h = \frac{1}{2} \left ( \int_{-\infty}^x
\rho_h(t,y) \ dy - \int_x^\infty \rho_h(t,y) \ dy \right ),$$
$$ \tilde{h} = \mathcal{F}[h],$$
$$ \tilde{g}_h = F - \tilde{h},$$
and
$$ \tilde{\rho}_h = \int \tilde{g}_h \ dv.$$

We estimate $\Vert \tilde{\rho}(t) - \tilde{\rho}_h(t) \Vert_p$ so
that we may apply a uniform bound to the iterates defined later.
As before, we consider $D$ defined by (\ref{D}) and find
$$ D ( \tilde{f} - \tilde{h} ) = (E - E_h) \partial_v \tilde{h}
= (E - E_h) (F^\prime - \partial_v \tilde{g}_h).$$ So, we have
\begin{equation}
\begin{array}{rcl}
\label{rhodiffeq} (\tilde{\rho} - \tilde{\rho}_h)(t,x) & = &
\int_0^t \int ( E - E_h) (s,X(s)) F^\prime(V(s)) \ dv \ ds \\
\\
& & - \int_0^t \int (E - E_h)(s,X(s)) \partial_v
\tilde{g}_h(s,X(s),V(s)) \ dv \ ds
\end{array}
\end{equation}
where $X(s)$ and $V(s)$ are defined by (\ref{char}).  Since $f,h \in
\mathcal{C}$, we use (\ref{E}) and (\ref{Ebound}) to find
$$\vert (E - E_h)(t,x) \vert \leq C \Vert \rho(t) - \rho_h(t)
\Vert_p$$ and
$$ \vert \partial_x(E - E_h)(t,x) \vert \leq C \Vert \rho(t) - \rho_h(t) \Vert_p
R^{-p}(x).$$ Using Lemma $2$, we have $$ \left \vert \int (E -
E_h) (s,X(s)) F^\prime(V(s)) \ dv \right \vert \leq C \Vert (\rho
- \rho_h)(s) \Vert_p R^{-p}(x).$$ Thus, it follows that
\begin{equation}
\label{EhF} \left \vert \int_0^t \int (E - E_h) (s,X(s))
F^\prime(V(s)) \ dv \ ds \right \vert \leq C \left( \int_0^t \Vert
\rho - \rho_h)(s) \Vert_p \ ds \right) R^{-p}(x).
\end{equation}
To estimate the remaining portion of (\ref{rhodiffeq}), we will
use :

\begin{lemma} For  $0 \leq s \leq t \leq T + \delta$ and $x \in
\mathbb{R}$,
$$\left \vert \int (E - E_h)(s, X(s)) \ \partial_v \tilde{g}_h(s,X(s),V(s)) \ dv \right \vert
\leq C \Vert (\rho - \rho_h)(s) \Vert_p R^{-p}(x).$$
\end{lemma}
It follows from Lemma $4$ that
\begin{equation}
\label{Ehg}
\left \vert \int_0^t \int (E - E_h)(s, X(s)) \
\partial_v \tilde{g}_h(s,X(s),V(s)) \ dv \ ds \right \vert \leq
C \left(\int_0^t \Vert (\rho - \rho_h)(s) \Vert_p \ ds \right)
R^{-p}(x).
\end{equation}
Using (\ref{EhF}) and (\ref{Ehg}) in (\ref{rhodiffeq}), we find
\begin{equation}
\label{rhohbound} \Vert (\tilde{\rho} - \tilde{\rho}_h)(t) \Vert_p
\leq C^{(3)} \left ( \int_0^t \Vert (\rho - \rho_h)(s) \Vert_p \
ds \right).
\end{equation}

Now, define the first iterate by

$$ f^{(0)}(t,x,v) = \left \{ \begin{array}{ll} f^{(T)}(t,x,v) & \mathrm{if} \ t \in [0,T] \\
f^{(T)}(T,x,v) & \mathrm{if} \ t \in [T,T+\delta]. \end{array}
\right.$$ Then, for $k = 0,1,2,3,...$, assuming $f^{(k)}$ is known,
define
\begin{eqnarray*}
g^{(k)} & = & F - f^{(k)}, \\
\rho^{(k)} & = & \int g^{(k)} \ dv, \\
E^{(k)} & = & \frac{1}{2} \left ( \int_{-\infty}^x \rho^{(k)}(t,y) \
dy - \int_x^\infty \rho^{(k)}(t,y) \ dy \right ), \\
f^{(k+1)} & = & \mathcal{F}[f^{(k)}].
\end{eqnarray*}
Since $f^{(0)}, f^{(1)} \in \mathcal{C}$, we have $$ \Vert (
\rho^{(1)} - \rho^{(0)} ) (t) \Vert_p \leq C^{(4)}I(t - T).$$
Then, we use (\ref{rhohbound}) to find $$ \Vert ( \rho^{(2)} -
\rho^{(1)} ) (t) \Vert_p \leq C^{(3)} C^{(4)} (t - T) I(t - T),$$
and
$$ \Vert ( \rho^{(3)} - \rho^{(2)} ) (t) \Vert_p \leq C^{(4)} \frac{1}{2} \left (C^{(3)} (t-T) \right )^2 I(t- T).$$ This
can be repeated and using induction, for any $k = 0,1,2,3,...$, we
have
$$\Vert ( \rho^{(k+1)} - \rho^{(k)} ) (t) \Vert_p \leq \frac{C^{(4)}
}{k!} \left ( C^{(3)} (t - T) \right )^k I(t - T).$$  Finally, for
$m,n \in \mathbb{Z}^+$ with $m > n$, we have
\begin{eqnarray*}
\Vert ( \rho^{(m)} - \rho^{(n)} ) (t) \Vert_p & \leq & \sum_{\ell
=
n}^{m-1} \Vert ( \rho^{(\ell + 1)} - \rho^{(\ell)} ) (t) \Vert_p \\
& \leq & \left ( \sum_{\ell = n}^{m-1} \frac{C^{(4)}}{\ell !} \left
( C^{(3)} (t - T) \right )^\ell \right ) I(t - T).
\end{eqnarray*}
Thus, $\rho^{(k)}$ is Cauchy with respect to the norm
$$ \Vert \sigma \Vert = \sup_{t \in [0, T + \delta]} \Vert \sigma(t) \Vert_p.$$
Using this estimate with (\ref{Ebound}), we may conclude that
$E^{(k)}$ is Cauchy in $L^\infty([0,T+\delta] \times \mathbb{R}).$
Define $X^{(k)}(s,t,x,v)$ and $V^{(k)}(s,t,x,v)$ by
\begin{equation} \label{chark} \left \{
\begin{array}{ll} \displaystyle \frac{dX^{(k)}}{ds} =
V^{(k)},  \ &
X^{(k)}(t,t,x,v) = x \\
\\
\displaystyle \frac{dV^{(k)}}{ds}, = - E^{(k)}
\vert_{(s,X^{(k)}(s))} \ & V^{(k)}(t,t,x,v) = v.
\end{array} \right.
\end{equation}
Then, since $f^{(n)}, f^{(m)} \in \mathcal{C}$, we use the Mean
Value Theorem to find
\begin{eqnarray*} \vert (X^{(n)} -
X^{(m)})(s) \vert + \vert (V^{(n)} - V^{(m)})(s)
\vert & \leq & \int_s^t ( \vert (V^{(n)} - V^{(m)})(\tau) \vert \\
& \ & \ \ + \vert
E^{(n)}(\tau,X^{(n)}(\tau)) - E^{(m)}(\tau,X^{(m)}(\tau)) \vert \ d\tau \\
& \leq & \int_s^t ( \vert (V^{(n)} - V^{(m)}(\tau) \vert  \\
& \ & \ \ + \vert E^{(n)}(\tau,X^{(n)}(\tau)) - E^{(n)}(\tau,X^{(m)}(\tau)) \vert \\
& \ & \ \ + \vert E^{(n)}(\tau,X^{(m)}(\tau)) -
E^{(m)}(\tau,X^{(m)}(\tau)) \vert) \ d\tau \\
& \leq & \int_s^t \left ( \vert (V^{(n)} - V^{(m)}(\tau) \vert
\right. \\
& \ & \ \ + \left. \Vert \rho^{(n)}(\tau) \Vert_p \ \vert (X^{(n)}
- X^{(m)}) (\tau) \vert \right
) \ d\tau \\
& \ & \ \ + \Vert E^{(n)} - E^{(m)} \Vert_{L^\infty([0,T+ \delta]
\times \mathbb{R})} \\
& \leq & \Vert E^{(n)} - E^{(m)} \Vert_{L^\infty([0,T+ \delta]
\times \mathbb{R})} \\
& \ & \ \ + C \int_s^t \left ( \vert (X^{(n)} - X^{(m)}) (\tau)
\vert + \vert (V^{(n)} - V^{(m)}(\tau) \vert \right ) \ d\tau
\end{eqnarray*}
Using Gronwall's Inequality, we have for $0 \leq s \leq t \leq T +
\delta$
$$ \vert (X^{(n)} - X^{(m)})(s) \vert  + \vert (V^{(n)} - V^{(m)})(s) \vert
\leq C \Vert E^{(n)} - E^{(m)} \Vert_{L^\infty([0,T+ \delta] \times
\mathbb{R})}.$$ Thus, $X^{(k)}$ and $V^{(k)}$ are uniformly Cauchy.
Then, by definition of the iterates,
$$ f^{(k+1)}(t,x,v) = f_0(X^{(k)}(0,t,x,v), V^{(k)}(0,t,x,v)).$$
So, define
\begin{equation}
\label{fdef} f(t,x,v)  = \lim_{k \rightarrow \infty}
f_0(X^{(k)}(0,t,x,v), V^{(k)}(0,t,x,v)).
\end{equation}
Then, since $f^{(k)} \in \mathcal{C}$ and thus $\Vert \rho^{(k)}(t)
\Vert_p \leq C$, we use (A-$3$), Lemma $1$, and the Bounded
Convergence Theorem to define
$$ \rho = \lim_{k \rightarrow \infty} \rho^{(k)} = \lim_{k
\rightarrow \infty} \int ( F - f^{(k)} ) \ dv = \int (F - f) \ dv.$$
Similarly, define $$ E = \lim_{k \rightarrow \infty} E^{(k)} =
\lim_{k \rightarrow \infty} \frac{1}{2} \left ( \int_{-\infty}^x
\rho^{(k)}(t,y) \ dy - \int_x^\infty \rho^{(k)}(t,y) \ dy \right )$$
and use Lebesgue's Dominated Convergence Theorem to conclude that
(\ref{E}) holds. Define
$$ X = \lim_{k
\rightarrow \infty} X^{(k)}$$ and
$$ V = \lim_{k
\rightarrow \infty} V^{(k)}.$$ It follows from (\ref{chark}) and the
field bound that (\ref{char}) holds.  Also, from (\ref{fdef}), we
have
\begin{equation}
\label{fdef2} f(t,x,v) = f_0(X(0,t,x,v), V(0,t,x,v))
\end{equation}
and thus
\begin{equation}
\label{fdef3}
f(t,x,v) = f(s,X(s,t,x,v),V(s,t,x,v))
\end{equation}
for all $s \in [0,t]$. Finally, $f^{(k)} \in C$ for each $k$ implies
$\rho$ is Lipschitz in $x$ and thus $E$ is $C^1$ in $x$. We may use
(\ref{ftildedef}) and Lemma $1$ to conclude that $E$ is continuous
in $t$. Thus, we see that $X$ and $V$ are $C^1$ from (\ref{char}),
$f$ is $C^1$ from (\ref{fdef2}), and $f$ satisfies the Vlasov
equation from
(\ref{fdef3}).\\

In order to show uniqueness, we let $f,h \in \mathcal{C}$ be given
solutions with $f(0,x,v) = h(0,x,v) = f_0(x,v)$ for all $x,v \in
\mathbb{R}$. Then, since $f,h \in \mathcal{C}$ solve (\ref{VP}) we
use (\ref{rhohbound}) to conclude
$$ \Vert (\rho - \rho_h)(t) \Vert_p \leq C (t - T) I(t-T).$$
We may repeatedly apply (\ref{rhohbound}) and use induction to find
for every $k \in \mathbb{N}$
$$ \Vert (\rho - \rho_h)(t) \Vert_p \leq \frac{C^{(4)}}{k!} \left ( C^{(3)} (t-T) \right )^k I(t-T).$$
This implies
$$ \Vert (\rho - \rho_h)(t) \Vert_p = 0.$$  Therefore, $\rho \equiv \rho_h$, $E \equiv E_h$
and by uniqueness of characteristics, $f \equiv h$.  Thus, the
proofs of Theorem $1$ and Theorem $2$ are complete.

\section*{Section 3}

The section which follows is devoted to the proof of Lemmas $1$,
$2$, $3$, and $4$.

\noindent \emph{\textbf{Proof of  Lemma $1$}} : Using (\ref{intE})
and (\ref{char}), we may conclude
\begin{eqnarray*}
\vert V(s,t,x,v) \vert & = & \left \vert v + \int_s^t
E(\tau,X(\tau,t,x,v)) d\tau \right \vert \\
& \geq & \vert v \vert - \int_0^t \Vert E(\tau) \Vert_\infty d\tau \\
& \geq & \vert v \vert - C^{(1)}.
\end{eqnarray*}
Similarly, we find $$ \vert V(s,t,x,v) \vert \leq \vert v \vert +
C^{(1)}.$$ Thus, for $\vert v \vert > 2 C^{(1)}$, we have
$$ \frac{1}{2} \vert v \vert \leq \vert v \vert - C^{(1)} \leq \vert V(s,t,x,v) \vert
\leq \vert v \vert + C^{(1)} \leq  \frac{3}{2} \vert v \vert$$ and
the proof of Lemma $1$ is
complete. \\

\vspace{0.25in}

\noindent \emph{\textbf{Proof of  Lemma $2$}} : Let $\mathcal{E}
\in \mathcal{C}^1(\mathbb{R})$, $H \in
\mathcal{C}^2_c(\mathbb{R})$, and $B
> 0$ be given with
\begin{equation}
\label{mcE}
\vert \mathcal{E}(x) \vert \leq B
\end{equation}
and
\begin{equation}
\label{mcEp}
\vert \mathcal{E}^\prime(x) \vert \leq B R^{-p}(x)
\end{equation}
for all $x \in \mathbb{R}$.  Let $t \in [0,T+ \delta]$ and $s \in
[0,t]$ be given.  Since $H$ has compact support, let $\tilde{W} :=
\sup \{ \vert v \vert : H(v) \neq 0 \}$ and define $C^{(5)} :=
2\max\{C^{(1)}, \tilde{W} \}.$  Then, using Lemma 1, we find for
$\vert v \vert
> C^{(5)}$ that $\vert V(s) \vert \geq \frac{1}{2} \vert v \vert
\geq \tilde{W}$, and thus $H^\prime(V(s)) = 0$.  Define
$$C^{(6)} := (T + \delta) (C^{(1)} + C^{(5)}).$$ We have for
$\vert x \vert \leq 2C^{(6)}$
\begin{eqnarray*}
\left \vert \int \mathcal{E}(X(s,t,x,v)) H^\prime(V(s,t,x,v)) \ dv
\right \vert & \leq & \int_{\vert v \vert \leq C^{(5)}} B \Vert
H^\prime \Vert_\infty \
dv \\
& \leq & CB \leq CBR^{-p}(x).
\end{eqnarray*}
Take $\vert x \vert > 2C^{(6)}$ and write

\begin{eqnarray*}
& \ & \int \mathcal{E}(X(s,t,x,v)) H^\prime(V(s,t,x,v))
\ dv \\
& = & \int_{\vert v \vert \leq C^{(5)}} \left [ \mathcal{E}(
X(s,t,x,v)) \left ( H^\prime(V(s,t,x,v)) -
H^\prime(v + \int_s^t E(\tau,x) \ d\tau) \right ) \right. \\
& \ & + ( \mathcal{E}(X(s,t,x,v)) - \mathcal{E}(x + (s-t)v) )
H^\prime(v + \int_s^t E(\tau,x) \ d\tau ) \\
& \ & + \frac{d}{dv} \left (H(v + \int_s^t E(\tau,x) \
d\tau) \ \mathcal{E}( x + (s-t)v) \right ) \\
& \ & - \left. H(v + \int_s^t E(\tau,x) \ d\tau) \ \frac{d}{dv}
\left ( \mathcal{E}(x + (s-t)v) \right ) \right ] \
dv \\
& =: & I + II + III + IV.
\end{eqnarray*}

To estimate $I$, we use the Mean Value Theorem to find $\xi_1$
between $x$ and $X(\tau)$ such that $$ E(\tau, X(\tau)) -
E(\tau,x) = \partial_x E(\tau, \xi_1) (X(\tau) - x).$$ In
addition, notice that for $\vert v \vert \leq C^{(5)}$, we have
using Lemma $1$
\begin{eqnarray*}
\vert \xi_1 \vert & \geq & \vert x \vert - \vert X(\tau) - x
\vert \\
& \geq & \vert x \vert - \left ( \int_\tau^t (\vert v \vert + C^{(1)}) \ ds \right )\\
& \geq & \vert x \vert - (T + \delta) (C^{(5)} + C^{(1)}) \\
& \geq & \vert x \vert - C^{(6)} \\
& \geq & \frac{1}{2} \vert x \vert.
\end{eqnarray*}
Using (\ref{rhobound}), (\ref{mcE}), and Lemma $1$ we find
\begin{eqnarray*}
I & \leq & CB \int_{\vert v \vert \leq C^{(5)}} \left \vert
H^\prime(V(s,t,x,v)) - H^\prime(v + \int_s^t E(\tau, x) \ d\tau) \right \vert \ dv \\
& \leq & CB \Vert H^{\prime \prime} \Vert_\infty \int_{\vert v
\vert \leq
C^{(5)}} \int_s^t \vert E(\tau, X(\tau)) - E(\tau,x) \vert \ d\tau \ dv \\
& \leq & CB\int_{\vert v \vert \leq C^{(5)}} \int_s^t \vert
\partial_x E(\tau, \xi_1) \vert \ \vert X(\tau) - x \vert \ d\tau \ dv \\
& \leq & CB \int_{\vert v \vert \leq C^{(5)}} \int_s^t \Vert
\rho(\tau) \Vert_p R^{-p}(\xi_1) \ \left ( \int_\tau^t \vert
V(\lambda) \vert \
d\lambda \ \right ) \ d\tau \ dv \\
& \leq & CB R^{-p} \left ( \frac{1}{2} \vert x \vert \right ) \left
( \int_{\vert v \vert \leq C^{(5)}} (\vert v \vert + C^{(1)}) \ dv
\right ) \left ( \int_s^t
(1 + r I(\tau-T)) \ d\tau \right ) \\
& \leq & CB R^{-p}(x).
\end{eqnarray*}
To estimate $II$, we again use the Mean Value Theorem and find
$\xi_2$ between $X(s)$ and $x + (s-t)v$ such that $$
\mathcal{E}(X(s)) - \mathcal{E}(x + (s-t)v) =
\mathcal{E}^\prime(\xi_2) (X(s) - (x + (s-t)v)).$$ In addition,
notice that for $\vert v \vert \leq C^{(5)}$, we use Lemma $1$ and
(\ref{intE}) to find
\begin{eqnarray*}
\vert \xi_2 \vert & \geq & \vert x + (s-t) v \vert - \vert X(s) -
(x + (s-t)v) \vert \\
& \geq & \vert x \vert - (T + \delta) C^{(5)} - \left \vert
\int_s^t \int_\tau^t E(\iota,X(\iota)) \ d\iota \ d\tau \right \vert \\
& \geq & \vert x \vert - (T + \delta)(C^{(5)} +
C^{(1)}) \\
& \geq & \frac{1}{2} \vert x \vert.
\end{eqnarray*}
Then, using (A-2), (\ref{intE}), and (\ref{mcEp}), we find
\begin{eqnarray*}
II & \leq &  C \Vert H^\prime \Vert_\infty \int_{\vert v \vert
\leq C^{(5)}} \vert \mathcal{E}(X(s,t,x,v)) -
\mathcal{E}(x + (s-t)v) \vert \ dv \\
& \leq & C \int_{\vert v \vert \leq C^{(5)}} \vert
\mathcal{E}^\prime (\xi_2) \vert \ \vert X(s) - x - (s-t)v \vert \
dv \\
& \leq & CB \int_{\vert v \vert \leq C^{(5)}} R^{-p}(\xi_2) \int_s^t
\int_\tau^t \vert E(\iota,X(\iota)) \vert \ d\iota \
d\tau \\
& \leq & CB R^{-p}(x).
\end{eqnarray*}
By the Fundamental Theorem of Calculus and compact support of $H$,
we find
$$ III = 0.$$ To estimate $IV$, notice that for $\vert v \vert
\leq C^{(5)}$,
\begin{equation}
\label{xstv} \vert x + (s-t)v \vert \geq \vert x \vert - (T +
\delta)C^{(5)} \geq \vert x \vert - C^{(6)} \geq \frac{1}{2} \vert
x \vert.
\end{equation}
Thus, we use (\ref{mcEp}), (\ref{xstv}), and Lemma $1$ to estimate
$IV$, which yields
\begin{eqnarray*}
IV & \leq & C \Vert H \Vert_\infty \int_{\vert v \vert
\leq C^{(5)}} \vert (s-t) \vert \ \vert \mathcal{E}^\prime(x+(s-t)v) \vert \ dv \\
& \leq & CB \int_{\vert v \vert \leq C^{(5)}} R^{-p}(x + (s-t)v) \ dv \\
& \leq & CBR^{-p}(x).
\end{eqnarray*}
Combining the estimates for $I$ - $IV$, we have $$ \left \vert
\int \mathcal{E}(X(s) H^\prime (V(s)) \ dv \right \vert \leq
CBR^{-p}(x)$$ for $\vert x \vert > 2C^{(6)}$, and the lemma
follows.

\vspace{0.25in}

\noindent \emph{\textbf{Proof of  Lemma $3$}} :  We must bound
$Q_{\tilde{g}}(t)$.  Define for $t \in [0,T + \delta]$,
$$Q(t) : = \sup \{ \vert v \vert : \exists \ x \in \mathbb{R},
\tau \in [0,t] \ s.t. \ \tilde{f}(\tau,x,v) \neq 0 \}.$$ By (A-1),
we know that $Q(0)$ is finite.  By the definition of $Q(t)$, if
$\vert V(0,t,x,v) \vert \geq Q(0)$, we have for every $y \in
\mathbb{R}$,
$$f_0(y,V(0,t,x,v)) = 0.$$ But, by Lemma $1$, if $\vert v \vert
\geq 2 \max\{Q(0), C^{(1)}\}$, then
$$ \vert V(0,t,x,v) \vert \geq \frac{1}{2} \vert v \vert \geq Q(0)$$
which implies that $f_0(y,V(0,t,x,v)) = 0$. So, if
$f_0(y,V(0,t,x,v)) \neq 0$ for some $y \in \mathbb{R}$, we must have
\begin{equation}
\label{vbound} \vert v \vert \leq 2 \max\{Q(0), C^{(1)}\}.
\end{equation}
The definition of $\tilde{f}$ implies that $\tilde{f}(t,x,v) =
f_0(X(0,t,x,v), V(0,t,x,v))$.  So, if $\tilde{f}(t,x,v) \neq 0$ for
some $t \in [0,T + \delta]$, $x,v \in \mathbb{R}$, then
(\ref{vbound}) must hold. Taking the supremum over $v$ of both sides
in (\ref{vbound}), we find
$$ Q(t) \leq 2 \max\{Q(0), C^{(1)}\} \leq C $$
for every $t \in [0,T + \delta]$.  Since $F$ is compactly supported
and $Q_{\tilde{g}}(t) = \max\{Q(t), W\}$, it follows that

\begin{equation}
\label{Qgtildebound} Q_{\tilde{g}}(t) \leq C
\end{equation} for all
$t \in [0,T + \delta]$, as well.

\vspace{0.25in}

\noindent \emph{\textbf{Proof of  Lemma $4$}} : Let $t \in
[0,T+\delta]$, $s \in [0,t]$, and $x \in \mathbb{R}$ be given.
Notice that by the definitions of $E$ and $E_h$, we have
$$\partial_x (E - E_h)(t,x) = (\rho - \rho_h)(t,x).$$
In addition, using (\ref{rhobound}) and (\ref{Ebound}), we have
\begin{equation}
\label{EmEh} \vert (E - E_h)(t,x) \vert \leq C \Vert (\rho -
\rho_h)(t) \Vert_p
\end{equation}
and
\begin{equation}
\label{rhomrhoh} \vert \partial_x (E - E_h)(t,x) \vert \leq C
\Vert (\rho - \rho_h)(t) \Vert_p R^{-p}(x).
\end{equation}
 Define $C^{(7)} := (T + \delta) (Q_{\tilde{g}_h}(T +
\delta) + C^{(1)})$. Then, for $\vert x \vert \leq 2C^{(7)}$, we use
(\ref{gtildetriple}) and (\ref{EmEh}) to find
\begin{eqnarray*}
\left \vert \int (E - E_h)(s, X(s)) \ \partial_v
\tilde{g}_h(s,X(s),V(s)) \ dv \right \vert & \leq & C \int_{\vert
v \vert \leq Q_{\tilde{g}_h}(s)} \Vert (\rho - \rho_h)(s) \Vert_p
\ dv \\
& \leq & C \Vert (\rho - \rho_h)(s) \Vert_p \\
& \leq & C \Vert (\rho - \rho_h)(s) \Vert_p R^{-p}(x).
\end{eqnarray*}

Now, let $\vert x \vert > 2 C^{(7)}$, and notice that for $\vert v
\vert \leq Q_{\tilde{g}_h}(s)$ and $\tau \in [0,s]$, we have
\begin{eqnarray*}
\vert X(\tau) \vert & \geq & \vert x \vert - \int_\tau^t \vert
V(\iota)
\vert \ d\iota \\
& \geq & \vert x \vert - \int_\tau^t \left ( \vert v \vert +
C^{(1)}
\right ) \ d\iota\\
& \geq & \vert x \vert - (T + \delta) (Q_{\tilde{g}_h}(T + \delta)
+ C^{(1)}) \\
& \geq & \vert x \vert - C^{(7)}.
\end{eqnarray*}
Thus, for $\vert v \vert \leq Q_{\tilde{g}_h}(s)$ and $\tau \in
[0,s]$, we have
\begin{equation}
\label{Xx} \vert  X(\tau) \vert \geq \frac{1}{2} \vert x \vert.
\end{equation}
Integrating along characteristics in (\ref{Ddvg}), we find
$$ \partial_v \tilde{g}_h \vert_{(s,X(s),V(s))} =
\partial_v \tilde{g}_h \vert_{(0,X(0),V(0))} - \left. \int_0^s
( E_h F^{\prime \prime} + \partial_x \tilde{g}_h ) \right
\vert_{(\tau, X(\tau), V(\tau))} \ d\tau.
$$ Using this, we may write
\begin{eqnarray*}
\left. \int [ (E - E_h) \partial_v \tilde{g}_h ] \right
\vert_{(s,X(s),V(s))} \ dv & = & \left. \int_{\vert v \vert \leq
Q_{\tilde{g}_h}(s)} (E - E_h) \right \vert_{(s,X(s),V(s))}
\\ & \ & \cdot \left [ \left.
\partial_v \tilde{g}_h \right \vert_{(0,X(0),V(0))} - \left. \int_0^s \left ( E_h
F^{\prime \prime} +
\partial_x \tilde{g}_h \right ) \right \vert_{(\tau, X(\tau), V(\tau))} \
d\tau \right ] dv \\
& \leq & I + II + III
\end{eqnarray*}
where
$$ I = \left \vert \int_{\vert v \vert \leq Q_{\tilde{g}_h}(s)} (E - E_h)(s,X(s)) \ \partial_v \tilde{g}_h (0,X(0),V(0)) \ dv \right \vert,$$
$$ II  = \left \vert \int_{\vert v \vert \leq Q_{\tilde{g}_h}(s)} (E - E_h)(s,X(s)) \int_0^s
E_h (\tau, X(\tau)) F^{\prime \prime}(V(\tau)) \ d\tau \ dv \right
\vert,$$ and
$$ III = \left \vert \int_{\vert v \vert \leq Q_{\tilde{g}_h}(s)} (E - E_h)(s,X(s)) \ \left (\int_0^s \partial_x \tilde{g}_h(\tau, X(\tau), V(\tau)) \
d\tau \right ) \ dv \right \vert .$$

To estimate $I$, we use (A-1), (\ref{EmEh}) and (\ref{Xx}) to find
\begin{eqnarray*}
I & \leq & C \Vert (\rho - \rho_h)(s) \Vert_p \int_{\vert v \vert \leq Q_{\tilde{g}_h}(s)} R^{-p}(X(0)) \ dv \\
& \leq & C \Vert (\rho - \rho_h)(s) \Vert_p Q_{\tilde{g}_h}(s) R^{-p}\left ( \frac{1}{2} \vert x \vert \right ) \\
& \leq & C \Vert (\rho - \rho_h)(s) \Vert_p R^{-p}(x) .
\end{eqnarray*}

To estimate $II$, we write
\begin{equation}
\label{addsubtract}
\begin{array}{rcl} (E - E_h)(s,X(s)) E_h (\tau, X(\tau)) F^{\prime
\prime}(V(\tau)) & = & \left ((E - E_h)(s,X(s)) - (E - E_h)(s,x)
\right ) E_h (\tau,
X(\tau)) F^{\prime \prime}(V(\tau)) \\
& \ & \ +  \ (E - E_h)(s,x) E_h (\tau, X(\tau)) F^{\prime
\prime}(V(\tau)).
\end{array}
\end{equation}
For the first term in (\ref{addsubtract}), we use the Mean Value
Theorem and find $\xi_3$ between $X(s)$ and $x$ such that $$ (E -
E_h)(s,X(s)) - (E - E_h)(s,x) = \partial_x(E - E_h)(s, \xi_3) \
(X(s) - x).$$ In addition, notice that for $\vert v \vert \leq
Q_{\tilde{g}_h}(s)$, we use Lemma $1$ and (\ref{intE}) to find
\begin{eqnarray*}
\vert \xi_3 \vert & \geq & \vert x \vert - \vert X(s) - x\vert \\
& \geq & \vert x \vert - \int_s^t \vert V(\iota) \vert \ d\iota \\
& \geq & \vert x \vert - (T + \delta) (Q_{\tilde{g}_h}(T + \delta) + C^{(1)})\\
& \geq & \vert x \vert - C^{(7)} \\
& \geq & \frac{1}{2} \vert x \vert.
\end{eqnarray*}
Thus, using (\ref{Ebound}) and (\ref{rhomrhoh}), we find
\begin{eqnarray*}
& & \left \vert \int_{\vert v \vert \leq Q_{\tilde{g}_h}(s)} ( (E -
E_h)(s,X(s)) - (E - E_h(s,x) ) E_h
(\tau, X(\tau)) F^{\prime \prime}(V(\tau)) \ dv \right \vert \ \ \ \ \\
& \leq & \int_{\vert v \vert \leq Q_{\tilde{g}_h}(s)} \left \vert
\partial_x(E - E_h)(s, \xi_3) \ (X(s) - x) \right \vert \ \vert
E_h (\tau, X(\tau)) F^{\prime \prime}(V(\tau)) \vert \ dv \\
& \leq & C (Q_{\tilde{g}_h}(s) + C^{(1)})(T + \delta) \ \Vert
F^{\prime \prime} \Vert_{L^\infty(\mathbb{R})} \Vert
\rho(\tau) \Vert_p  \Vert (\rho - \rho_h)(s) \Vert_p R^{-p}(\xi_3) \\
& \leq & C \Vert \rho(\tau) \Vert_p \Vert (\rho - \rho_h)(s)
\Vert_p R^{-p}(x)
\end{eqnarray*}
and therefore,
\begin{equation}
\label{firstterm} \left \vert \int_{\vert v \vert \leq
Q_{\tilde{g}_h}(s)} ( (E - E_h)(s,X(s)) - (E - E_h(s,x) ) E_h (\tau,
X(\tau)) F^{\prime \prime}(V(\tau)) \ dv \right \vert \leq C \Vert
\rho(\tau) \Vert_p \Vert (\rho - \rho_h)(s) \Vert_p R^{-p}(x).
\end{equation}
Then, to estimate the second term of (\ref{addsubtract}) we use
(\ref{Ebound}), (\ref{EmEh}), and (\ref{rhomrhoh}) in Lemma $2$
(with $\mathcal{E} = (E - E_h)(s,x) E_h(\tau, X(\tau))$ and $H =
F^\prime$) to find

\begin{equation}
\label{secondterm} \left \vert \int_{\vert v \vert \leq
Q_{\tilde{g}_h}(s)} (E - E_h)(s,x) E_h (\tau, X(\tau)) F^{\prime
\prime}(V(\tau)) \ dv \right \vert \leq C \Vert \rho_h (\tau)
\Vert_p \ \Vert (\rho - \rho_h)(s) \Vert_p \ R^{-p}(x).
\end{equation}
Finally, using (\ref{rhobound}), (\ref{firstterm}), and
(\ref{secondterm}), we find
\begin{eqnarray*}
II & = & \left \vert \int_0^s \int_{\vert v \vert \leq
Q_{\tilde{g}_h}(s)} (E - E_h)(s,X(s)) E_h (\tau, X(\tau))
F^{\prime \prime}(V(\tau)) \ dv \ d\tau \right \vert \\
& \leq & \int_0^s \left ( \left \vert \int_{\vert v \vert \leq
Q_{\tilde{g}_h}(s)} ((E - E_h)(s,X(s)) - (E - E_h)(s,x)) E_h
(\tau, X(\tau)) F^{\prime \prime}(V(\tau)) \ dv \right \vert \right. \\
& \ &  +  \ \left. \left \vert \int_{\vert v \vert \leq
Q_{\tilde{g}_h}(s)} (E - E_h)(s,x) E_h (\tau, X(\tau))
F^{\prime \prime}(V(\tau)) \ dv \right \vert \right ) d\tau \\
& \leq & C \Vert (\rho - \rho_h)(s) \Vert_p \left ( \int_0^s \Vert
\rho_h (\tau) \Vert_p  \ d\tau \right ) \ R^{-p}(x) \\
& \leq & C \Vert (\rho - \rho_h)(s) \Vert_p R^{-p}(x) \left (
\int_0^s (1 + r I(\tau-T)) \ d\tau \right ) \\
& \leq & C \Vert (\rho - \rho_h)(s) \Vert_p R^{-p}(x).
\end{eqnarray*}

To estimate $III$, we use (\ref{gtildetriple}), (\ref{EmEh}), and
(\ref{Xx}) to find
\begin{eqnarray*}
III & \leq & C \int_{\vert v \vert \leq Q_{\tilde{g}_h}(s)} \Vert
(\rho - \rho_h)(s) \Vert_p \left ( \int_0^s \Vert \partial_x
\tilde{g}_h(\tau) \Vert_p R^{-p}(X(\tau)) \ d\tau \right ) \ dv \\
& \leq & C Q_{\tilde{g}_h}(s) \ \left ( \int_0^s \Vert
\partial_x \tilde{g}_h(\tau) \Vert_p \ d\tau \right ) \  \Vert (\rho - \rho_h)(s)
\Vert_p R^{-p} \left( \frac{1}{2} \vert x \vert \right ) \\
& \leq & C \Vert (\rho - \rho_h)(s) \Vert_p R^{-p}(x)
\end{eqnarray*}
Combining the estimates for $I$, $II$, and $III$, the lemma
follows and the proof is complete.

\end{document}